\documentclass{article}

\title{A note on multiple Seshadri constants on surfaces.}
\author{Luis Fuentes Garc\'{\i}a}
\date{}

\newtheorem{teo}{Theorem}[section]

\newtheorem{prop}[teo]{Proposition}
\newtheorem{cor}[teo]{Corollary}
\newtheorem{lemma}[teo]{Lemma}

\def\g2{\pi}

\def\CC{\leavevmode\hbox{$\rm I\!\!\!C$}}

\def\impp{{\quad\Rightarrow\quad}}

\def\qed{\hspace{\fill}$\rule{2mm}{2mm}$}

\begin{document}

\maketitle

\vspace{0.1cm}

\begin{abstract}

We give a bound for the multiple Seshadri constants on surfaces with Picard number $1$. The result is a natural extension of the bound of A. Steffens for simple Seshadri constants. In particular, we prove that the Seshadri constant $\epsilon(L; r)$ is maximal when $rL^2$ is a square.

{\bf MSC (2000):} Primary 14C20; secondary, 14J60.

 {\bf Key Words:} Seshadri constants, algebraic surfaces.

\end{abstract}

\section{Introduction.}

The multiple Seshadri constants are a natural generalization of the Seshadri constants at single points defined by Demailly in \cite{De92}. If $X$ is a smooth projective surface, $L$ is a nef bundle on $X$ and $P_1,\ldots,P_r$ are distinct points in $X$, then the Seshadri constant of $L$ at $P_1,\ldots,P_r$ is:
$$
\epsilon(L; P_1,\ldots,P_r)=inf \frac{C\cdot L}{\sum_{i=1}^r mult_{P_i}C}
$$
where $C$ runs over all curves passing through at least one of the points $P_1,\ldots,P_r$. When the points are general we will write $\epsilon(L; r)$. This constants have the upper bound:
$$
\epsilon(L; r)\leq \sqrt{\frac{L^2}{r}}.
$$
However, explicit values are difficult to compute even when $r=1$. In \cite{St98}, A. Steffens proved the following result for simple Seshadri constants on surfaces with Picard number $1$:

\begin{prop}

Let $X$ be a surface with  $\rho(X)=rk(NS(X))=1$ and let $L$ be an ample generator of $NS(X)$. Let $\alpha$ be an integer with $\alpha^2\leq L^2$. If $x\in X$ is a very general point, then $\epsilon(L,x)\leq \alpha$. In particular, if $\sqrt{L^2}$ is an integer, then $\epsilon(L,x)=\sqrt{L^2}$.

\end{prop}

Some results in the same direction have been proved for multiple Seshadri constants. In \cite{Ha03}, Harbourne defines:
$$
\varepsilon_{r,l}=max\left\{\frac{\lfloor  d\sqrt{rl} \rfloor}{dr}, \quad 1\leq d\leq \sqrt{\frac{r}{l}}\right \}\cup
\left\{\frac{1}{\left\lceil\sqrt{\frac{r}{l}}\right\rceil}\right\}\cup
\left\{\frac{dl}{\lceil d\sqrt{rl}\rceil}, \quad 1\leq d\leq \sqrt{\frac{r}{l}}\right \}
$$
and he shows the following bound:

\begin{teo}

Let $l=L^2$, where $L$ is a very ample divisor on an algebraic surface $X$. Then $\epsilon(L; r)\geq \varepsilon_{l,r}$, unless $l\leq r$ and $rl$ is a square, in which case $\sqrt{l,r}=\varepsilon_{r,l}$ and $\epsilon(L; r)= \sqrt{l/r}$.

\end{teo}

When $l\leq r$ and $L$ is very ample this implies:
$$
\epsilon(L; r)\geq \frac{[\sqrt{rL^2}]}{r}.
$$
Moreover, if $rL^2$ is a square, $\epsilon(l; r)$ is maximal.

On the other hand, J. Roe in \cite{Ro04} relates the simple and multiple Seshadri constants. As a consequence of his main theorem and the result of Steffens, he obtains:

\begin{cor}

Let $X$ be a smooth projective surface defined over $\CC$,  $L$ an ample generator of $NS(X)$ and $r\geq 9$. The Nagata's conjecture implies:
$$
\epsilon(L,r)\geq \frac{[\sqrt{L^2}]}{\sqrt{r}}.
$$

\end{cor}

In this note, we extends the result of Steffens for multiple Seshadri constants. We prove:

\begin{teo}

Let $X$ be a surface with  $\rho(X)=rk(NS(X))=1$ and let $L$ be an ample generator of $NS(X)$. Then
$$
\epsilon(L; r)\geq \frac{[\sqrt{rL^2}]}{r}.
$$
In particular, if $rL^2$ is a square, $\epsilon(L;r)$ is maximal.

\end{teo}

In this case, the Harbourne's hypothesis of very ampleness of $L$ and $L^2\leq r$ are not necessary. Furthermore, we do not need to use the Nagata's conjecture. 

The proof of the theorem is a natural generalization of the method of Steffens.

\section{Proof of the Theorem.}

Let $X$ be a surface with $\rho(X)=1$ and let $L$ be an ample generator of $NS(X)$. Let  $\alpha$ be an integer  with $\alpha^2<rL^2$. Let us suppose that:
$$
\epsilon(L;r)<\frac{\alpha}{r}\leq \sqrt{\frac{L^2}{r}}.
$$
Then, there is a Seshadri exceptional curve $C$ with multiplicities $(m_1,\ldots,m_r)$ at very general points, such that:
$$
\epsilon(L,r)=\frac{C\cdot L}{M}, \quad \hbox{where} \quad M=\sum_{i=1}^{r} m_i.
$$
In order to bound this multiplicities, let us recall two lemmas.

\begin{lemma}

Let $X$ be a smooth surface  and let $(C_t,(P_1)_t,\ldots,(P_r)_t)$ be a one parameter
family of pointed curves on $X$ with $mult_{(P_i)_t}(C_t)=m_i$. Then:
$$
C_t^2\geq \sum_{i=1}^r m_i^2-min(m_1,\ldots,m_r).
$$

\end{lemma}
{\bf Proof:} See \cite{Xu95}. \qed

\begin{lemma}\label{almost}

Let $(X,L)$ be a polarized surface with Picard number $1$ and let $P_1,\ldots,P_r$ be general points on $X$. If $\epsilon(L; P_1,\ldots,P_r)<\sqrt{L^2/r}$ then any irreducible Seshadri curve is almost-homogeneous.

\end{lemma}
{\bf Proof:} See \cite{StSz04}. \qed

\begin{cor}

With the previous notation:
$$
rC^2\geq M(M-1).
$$

\end{cor}
{\bf Proof:} Applying the two lemmas, we know that  the multiplicities of $C$ are $(m_1,\ldots,m_r)=(a,\ldots,a,b)$ and
$$
rC^2\geq r(r-1)a^2+rb^2-r \, min(a,b).
$$
From this:
$$
\begin{array}{l}
{rC^2- M(M-1)\geq}\\
{\qquad  \geq r(r-1)a^2+rb^2-r \, min(a,b)-((r-1)a+b)^2+(r-1)a+b=}\\
{\qquad =(r-1)a^2+(r-1)b^2-2ab(r-1)+(r-1)a+b-r\,  min(a,b)=}\\
{\qquad =(r-1)(a-b)^2+(r-1)a+b-r \, min(a,b).}\\
\end{array}
$$
If $a\geq b$ then it holds:
$$
rC^2-M(M-1)\geq (r-1)((a-b)^2-(a-b))\geq 0.
$$
When $a<b$:
$$
rC^2-M(M-1)\geq (r-1)(a-b)^2+b-a\geq 0.
$$ \qed

Now, we can extend the proof of Steffens.  We have:
$$
\epsilon(L;r)<\frac{\alpha}{r} \impp \frac{C\cdot L}{M}<\frac{\alpha}{r} \impp rC\cdot L<\alpha M.
$$
On the other hand, since  $\rho(X)=1$, there is an integer $d$ such that $C\equiv dL$ and:
$$
rdL^2<\alpha M\impp \alpha^2d<\alpha M\impp \alpha d<M \impp \alpha d\leq M-1.
$$
Thus, applying the bound of the previous Corollary:
$$
M(M-1)\leq rC^2=rdC\cdot L<\alpha d M\leq M(M-1),
$$
and this is a contradiction. \qed

{\bf E-mail:} lfuentes@udc.es

Luis Fuentes Garc\'{\i}a. 

Departamento de M\'etodos Matem\'aticos y Representaci\'on.

E.T.S. de Ingenieros de Caminos, Canales y Puertos. 

Universidad de A Coru\~na. Campus de Elvi\~na. 15192 A Coru\~na (SPAIN)


\begin{thebibliography}{77}

\bibitem{De92}{\sc Demailly, J.-P.}
{\it Singular hermitian metrics on positive line bundles}
Lecture Notes Math. {\bf 1507}, 87-104 (1992).

\bibitem{Ha03}{\sc Harbourne, B.}
{\it Seshadri constants and very ample divisors on algebraic surfaces.}
J. Reine Angew. Math {\bf 559}, 115-122 (2003).	

\bibitem{Ro04}{\sc Roe, J.}
{\it A relation between one-point and multi-point Seshadri constants.}
J. Algebra {\bf 274}, Nº2, 643-651 (2004).

\bibitem{St98}{\sc Steffens, A.} 
{\it Remarks on Seshadri constants.}
Math. Z. {\bf 227}, 505-510 (1998).

\bibitem{StSz04}{\sc Strycharz-Szemberg, B.; Szemberg, T.}
{\it Remarks on the Nagata conjecture.}
Serdica Math. J. {\bf 30}, Nº 2-3, 405-430 (2004).

\bibitem{Xu95}{\sc Xu, G.}
{\it Ample line bundles on smooth surfaces.}
J. Reine Angew. Math. {\bf 469}, 199-209 (1995).



\end{thebibliography}
\end{document}